\documentclass[12pt,draft]{article}
\usepackage{a4wide}
\usepackage{graphicx}
\usepackage{amsfonts}
\usepackage{pdfsync}
\usepackage{amsmath}
\usepackage{enumerate}
\usepackage{amsthm}
\usepackage{float}

\usepackage{pgf,tikz,pgfplots}
\pgfplotsset{compat=1.15}
\usepackage{mathrsfs}
\usetikzlibrary{arrows}

\usepackage[colorlinks=true,linkcolor=blue,citecolor=blue]{hyperref}

\newcommand{\R}{\mathbb{R}}

\newtheorem{theorem}{Theorem}[section]
\newtheorem{lemma}{Lemma}[section]
\newtheorem{proposition}{Proposition}[section]
\newtheorem{corollary}{Corollary}[section]
\newtheorem{definition}{Definition}[section]
\newtheorem{remark}{Remark}[section]

\newcounter{exam}
\renewcommand{\theexam}{\thesection.\arabic{exam}}
\newenvironment{example}{
	\refstepcounter{exam}\paragraph{Example \theexam.\hspace*{-0.6em}}
	}
{$\Box$}

\title{Translators of the Mean Curvature Flow in Hyperbolic Einstein's Static Universe}

\author{Miguel Ortega and Buse Yalç\i n}
\date{\today}

\begin{document}

\maketitle

\begin{abstract}
In this study, we deal with non-degenerate translators of the mean curvature flow in the well-known hyperbolic Einstein's static universe. We classify translators foliated by horospheres and rotationally invariant ones, both space-like and time-like. For space-like translators, we show a uniqueness theorem as well as a result to extend an isometry of the boundary of the domain to the whole translator, under simple conditions. As an application, we obtain a characterization of the the bowl when the boundary is a ball, and of certain translators foliated by horospheres whose boundary is a rectangle.
\end{abstract}


\maketitle

\section{Introduction}\label{sec1}
The dynamic interplay between mean curvature flow (MCF) and the background of the hyperbolic Einstein's static universe (HESU) forms the core of our exploration in this study. The MCF, as a geometric evolution process, has been a pivotal subject of study with applications ranging from physics to mathematics. The study of the MCF traditionally focus on hypersurfaces within Euclidean space $\mathbb{R}^{n}$.  One approach involves examining specific solutions know as translating solutions or translators, which remain invariant under a subgroup of translations in the ambient space. The fundamental tool is to simplify the MCF to the equation
\begin{equation}\label{originalPDE}
\textbf{H}=v^{\perp},
\end{equation} 
where $n\geq2$, \textbf{H} represents the mean curvature vector of the immersion, $v$ is a constant (unit) vector, and $v^{\perp}$ denotes the orthogonal projection of $v$ into the normal bundle of the hypersurface. Some authors have simplified the equation \eqref{originalPDE} by examining rotationally invariant hypersurfaces, see \cite{CSS}. Moreover, in \cite{LO} they have adopted a broader approach, incorporating a cohomogeneity one action on $M$. This action involves the isometries of a Lie group, and it is structured such that the resulting hypersurfaces have constant mean curvature, with at most two exceptions (refer to \cite{AA} for further information). For example, in \cite{Bueno1} and \cite{Bueno2}, Bueno considered rotationally invariant translators in the product of the real hyperbolic plane and a real line, denoted as $\mathbb{H}^2\times\R$. In \cite{LM}, Lira and Mart\'in generalized the study by Bueno to Riemannian products $M\times\mathbb{R}$. Also, Pipoli considered translators in the solvable group $\mathbf{Sol}_3$ in \cite{Sol} and the Heisenberg group in \cite{Heis}. On the other hand, Kim in \cite{coreano} moved to the Minkowski space, again using some groups of isometries on space-like translators. This Lorentzian setting was studied in much more general Lorentzian spaces in \cite{LO}. Some other similar problems have been extensively examined in the existing literature (\cite{AW},\cite{CSS},\cite{HIMW},\cite{MSHS15},..., and references therein).

Our approach involves translators of the MCF in HESU. This study takes into account both the hyperbolic nature and static structure of the HESU. The involved techniques in our study are basically three, namely, the action of Lie groups of isometries to simplify some PDE into an ODE, the use of dynamical systems to solve a boundary problem, and general theory of elliptic quasilinear PDE. That is to say, this provides a framework for understanding translations of MCF in the context of the HESU, aiming to classify specific instances that emerge in this universe. 

The structure of this paper is the following: In Section 2, we review basic tools, highlighting that a graphical translator is characterized by a function satisfying the PDE \eqref{PDE 2}. We also exhibit a big family of time-like translators in HESU which can be constructed from minimal or totally geodesic hypersurfaces in the hyperbolic space.

In Section 3, we present the classification and properties of space-like and time-like translators foliated by horospheres, in Theorems \ref*{Theorem 3.1} and \ref*{Theorem 3.2}.

Section 4 is devoted to exploring rotationally invariant translators, which are those invariant by $SO(n)$. To facilitate it, we first recall the Lorentz space $\mathbb{L}^{n+1}$ as an essential mathematical background. We obtain a full classification of these translators, both space-like and time-like, obtaining 5 types in total, in Theorems \ref{theorem 4.1} and \ref{4.2}. An important space-like example  is known as the \textit{bowl}. Except the bowl, all of them exhibit one conic singularity. The only time-like example is a bigraph, known as the \textit{spindle}, and defined through the construction outlined in Lemma \ref{lemma 4.12}.

Finally, in Section 5, as an application of our findings, we show a uniqueness result for space-like translators based on the boundary, in Theorem \ref{uniq}. The methods employed throughout involves tools from quasilinear elliptic partial differential equations and isometries. As an consequence, we  prove that an isometry of the sphere can be extended to the translator under simple conditions. With this tool, we prove that a space-like translator whose boundary is a round ball has to be a compact piece of a bowl. In our last result, we characterize certain translators foliated by horospheres whose boundary is a rectangle in Corollary \ref{Corollary 5.4}.

\section{Preliminaries}\label{setup}
Let $M=\mathbb{H}^{n}$, $n\geq2$, be the hyperbolic space with its usual hyperbolic metric $g_{\mathbb{H}^{n}}$. We define the metric $\langle,\rangle$ as the usual metric in Euclidean space. We consider the Hyperbolic Einstein Static Universe (HESU), namely, the product $\bar{M}=\mathbb{H}^{n}\times\mathbb{R}$ with Lorentzian metric $g=g_{\mathbb{H}^{n}}-dt^{2}$. Take $\left(p,t\right)\in\mathbb{H}^{n}\times\mathbb{R}$. Given an open subset $\Omega$ of $\mathbb{H}^{n}$, we consider a function $u\in C^{\infty}\left(\Omega,\mathbb{R}\right)$ and construct its graph map $\Gamma:\Omega\rightarrow\bar{M}$, where $\Gamma(p)=\left(p,u(p)\right)$. Given the metric $\gamma=\Gamma^{*}\langle,\rangle$ on $\Omega$, we assume that $F:\left(\Omega,\gamma\right)\rightarrow\left(\bar{M},g\right)$ is a non-degenerate hypersurface. Under the usual identifications, for each $X\in TM$, we have 
\begin{equation*}
	dF(X)=\left(X,du(X)\right)=\left(X,\langle\nabla u,X\rangle\right),
\end{equation*} where $\nabla u$ is the $\langle,\rangle$-gradient of $u$. The upward vector field is 
\begin{equation}
	\nu=\frac{1}{W}\left(\nabla u,1\right),\qquad W=+\sqrt{\varepsilon\left(\lvert\nabla u \rvert^{2}-1\right)},
\end{equation} where $\varepsilon :=\mathrm{sign}\left(\lvert\nabla u \rvert^{2}-1\right)= \pm1$ is a constant function on the whole $\Omega$. Note that $g\left(\nu,\nu\right)=\varepsilon$. The following proposition is known (see for example \cite{LO}).
\begin{proposition}\label{proposition 2.1}
Under the previous setting, $\Gamma$ is a graphical translator if, and only if, function u satisfies the quasilinear PDE
	\begin{equation}\label{PDE 2}
		\mathrm{div}\left(\frac{\nabla u}{\sqrt{\varepsilon\left(\lvert\nabla u \rvert^{2}-1\right)}}\right)=\frac{1}{\sqrt{\varepsilon\left(\lvert\nabla u \rvert^{2}-1\right)}} .
	\end{equation}
\end{proposition}
\begin{example}
	Take $M$ a minimal or a totally geodesic hypersurface in $\mathbb{H}^{n}$. Then, the product $M\times\mathbb{R}$ is also a time-like translator in $\mathbb{H}^{n}\times_{-1}\mathbb{R}$, and it is not graphical. Indeed, since it is minimal or a totally geodesic in $\mathbb{H}^{n}$, then $M\times\mathbb{R}$ is also minimal or a totally geodesic, so its mean curvature vector $H=0$. But now, $\partial_{t}\in T\left(M\times\mathbb{R}\right)$. Then, $H=0=\left(\partial t\right)^{\perp}$, the projection of $\partial_{t}$ onto $T\left(M\times\mathbb{R}\right)$.                                                                
\end{example}
\section{Translators foliated by Horospheres}\label{horospheres}
 Before going any further, let us introduce the following space $\left(\mathbb{R}^{k},\langle, \rangle \right)=$ \textit{Flat Euclidean Space}.
We use the following model of the hyperbolic space

\begin{equation*}
	\mathbb{H}^{n}=\left\{\left(x_1,\ldots,x_n\right)\in \mathbb{R}^{n}\vert x_1>0 \right\},\, n\geq2,
\end{equation*}
equipped with the hyperbolic metric $g_{\mathbb{H}^{n}}=g_p=\frac{1}{{x_1}^{2}} \langle,\rangle_p$, $p\in\mathbb{H}^{n}$. The $(n-1)$-dimensional Horosphere  in $\mathbb{H}^{n}$ is given by
\begin{equation*}
	H^{n-1}=\left\{p=\left(x_1,\ldots,x_n\right)\in \mathbb{H}^{n}\vert x_1=\hat{x}_1\right\},\, \hat{x}_1>0.
\end{equation*}
We give a normal vector $N: H^{n-1} \longrightarrow T\mathbb{H}^{n},\, p\in\mathbb{H}^{n-1}\mapsto  x_{1}\partial_{1}\vert_p$. 
The tangent space $T_p H^{n-1}$ of $H^{n-1}$ can be written as
\begin{equation*}
	p\in H^{n-1},\, T_p H^{n-1}=\left\{X\in \mathbb{R}^{n}\vert X_1=0 \right\},
\end{equation*}
where for all $X\in T_p H^{n-1}$, one has $g\left(X,N\right)=0$.
Let $\nabla$ be the Levi-Civita connection on $\mathbb{H}^{n}$ and $A$ is the Shape operator of $N$ (Weingarten's operator). Given $X=\sum_{i=2}^{n}X_{i}\partial_{i}$ $\left(\mathrm{because} \, X_{1}=0, \, X\in TH^{n-1}\right)$, a simple computation gives $$AX=X.$$
We get $p\in H^{n-1},\, h(p)=\mathrm{trace}_g(A)=\left(n-1\right)$.

Take $\Omega\subset\mathbb{H}^{n}$ open and connected, $u\in C^{\infty}\left(\Omega,\mathbb{R}\right)$. Let $\Gamma:\Omega=\mathring{\Omega}\subset\mathbb{H}^{n} \longrightarrow M,\, p\mapsto\left(p,u(p)\right)$ be a translator, as in Section \ref{setup}. We wish the graph of $\Gamma$ to be foliated by horospheres. Therefore, $\Omega$ has to be also foliated by horospheres, that is to say $\Omega=I\times\mathbb{R}^{n-1}$, $I\subset\mathbb{R}$ an interval. We wish $u$ to be invariant by horospheres, so that $u$ only depends on one variable, namely,
\begin{equation*}
	u:\Omega\rightarrow\mathbb{R},\quad u:I\times\mathbb{R}^{n-1}\rightarrow\mathbb{R},\,\left(x_{1},\ldots ,x_{n}\right)\mapsto u(x_1).
\end{equation*}
Let $\Sigma$ be the Lie group such that the orbits of $\Sigma$ in $\mathbb{H}^{n}$ are just the horospheres, that is to say
$\Sigma=\left(\mathbb{R}^{n-1},+\right)$. The map
\begin{equation*}
	\begin{aligned}
		\Phi:\mathbb{R}^{n-1}\times\mathbb{H}^{n}&\rightarrow\mathbb{H}^{n},\\
		\left(\left(v_{1},\ldots,x_{n-1}\right),\left(x_{1},\ldots,v_{n}\right)\right)&\mapsto\left(x_{1},x_{2}+v_{1},\ldots,x_{n}+v_{n-1}\right),
	\end{aligned}
\end{equation*}
is an smooth action of $\Sigma$ on $\mathbb{H}$ by isometries. The projection map is
\begin{equation*}
	\rho:\mathbb{H}^{n}\rightarrow\mathbb{R},\,\rho\left(x_{1},\ldots,x_{n}\right)=\ln\left(x_{1}\right).
\end{equation*}
We put $\rho\left(\Omega\right)=I$. In particular,
$f:I\subset\mathbb{R}\rightarrow\mathbb{R},$ $u=f\circ\rho$ such that $F\left(p\right)=\left(p,f(\rho(p))\right)$, and then $\mathrm{graph}\left(f\circ\rho\right)$ is foliated by horospheres. Take $s\in\mathbb{R},\, H_{s}=\rho^{-1}(s)$ one horosphere. We compute the mean curvature of $H_{s}$ with respect to $-\nabla\rho$.
\begin{equation*}
h=\mathrm{div}\left(\nabla\rho\right)= \sum_{i=1}^{n-1} g\left(\nabla_{e_{i}}\nabla\rho,e_{i}\right)=n-1,
\end{equation*} where $\left(e_{1},\ldots,e_{n-1},e_{n}=\nabla\rho\right)$ is a local orthonormal frame on $T\mathbb{H}^{n}$ on each level set. According to \textit{Theorem 3.5} of paper \cite{LO}, $\Gamma_{u}$ is a translator if, and only if,
\begin{equation}
	f''\left(s\right)=\left(1-f'(s)^{2}\right)\left(1-f'(s)(n-1)\right),\quad n\geq2.
\end{equation}
\textbf{Case I.}
If $\varepsilon=-1$ then, $0<(-1)\left(-1+f'(s)^{2}\right)$ therefore, $\lvert f'(s)\rvert<1$, $\Gamma$ is space-like.
\newline
\textbf{Case II.}
If $\varepsilon=+1$ then, $0<(+1)\left(-1+f'(s)^{2}\right)$ therefore, $\lvert f'(s)\rvert>1$, $\Gamma$ is time-like.

Now, let us take z, we get

\begin{equation}
	w'(s)=(1-w(s)^2)(1-mw(s))\quad\mathrm{where}\quad m \in \mathbb R, \quad m\geq1,
\end{equation} From now, we will discuss all possible solutions to this ODE.
\subsection{Case $w(s)=\pm1$}
Function $f$ becomes $f(s)=\pm s+f_{0}$. We discard this case, because the associated hypersurfaces $\Gamma$ are degenerate.

\subsection{Case $w(s)=\frac{1}{m}$}
We obtain the function $f_{1}(s)=\frac{s}{m}+f_{0}$. 

\subsection{Case $w(s)\neq\pm1$ and $w(s)\neq\frac{1}{m}$}
We make some computations:
\begin{equation*}
	\frac{w'(s)}{(1-w(s)^2)(1-mw(s))} = 1  \Rightarrow \int \frac{w'(s)}{(1-w(s)^2)(1-mw(s))}ds =s-s_0
\end{equation*}
for some $s_0 \in \mathbb{R}$. Take the change of variable $x=w(s)$,
\begin{equation}\label{5}
	\int\frac{dx}{(1-x^2)(1-mx)}=s-s_0 .
\end{equation}
\subsubsection{Assume that m=1}
We have from \eqref{5}
\begin{equation}
\mathbf{P(x)}:=\int\frac{dx}{(1-x^2)(1-x)}=- \frac{1}{4} \ln\lvert1-x\lvert+\frac{1}{2}\cdot\frac{1}{1-x}+\frac{1}{4}\ln\lvert1+x\lvert.
\end{equation}
Clearly, we need the following intervals
$\mathbb R- \left\{-1, +1\right\}=\left(-\infty,-1\right)\cup\left(-1,+1\right)\cup\left(+1,+\infty\right)$.
 and we obtain different functions, by restricting function $\mathbf{P}$ to the different intervals, namely, $P_1=\mathbf{P}\vert_{\left(-\infty,-1\right)}$, $P_2=\mathbf{P}\vert_{\left(-1,+1\right)}$, $P_3=\mathbf{P}\vert_{\left(+1,+\infty\right)}$. We show some properties of them:
\begin{enumerate}
	\item $\lim_{x\rightarrow-\infty}P_1(x)=0$, $\lim_{\begin{subarray}{l}
			x\rightarrow-1\\
			x<-1
	\end{subarray}}P_1(x)=-\infty$; $P_{1}$ strictly decreasing, $\mathrm{Im}(P_1)=(-\infty,0)$.
\item $\lim_{\begin{subarray}{l}
		x\rightarrow-1\\
		x>-1
\end{subarray}}P_2(x)=-\infty$, $\lim_{\begin{subarray}{l}
		x\rightarrow1\\
		x<1
\end{subarray}}P_2(x)=+\infty$; $P_{2}$ strictly increasing, $\mathrm{Im}(P_2)=\mathbb{R}$.
\item $\lim_{\begin{subarray}{l}
		x\rightarrow1\\
		x>1
\end{subarray}}P_3(x)=-\infty$, $\lim_{+\infty}P_3(x)=0$; $P_{3}$ strictly increasing, $\mathrm{Im}(P_3)=(-\infty,0)$.
\end{enumerate}

\begin{figure}[h]
	\begin{minipage}[t]{0.49\linewidth}
		\centering
		\begin{tikzpicture}[scale = 1]
			\draw (-5,0)--(0.5,0);
			\draw (0,0.5) -- (0,-3);
			\draw[dashed] (-1,0.2) -- (-1,-3);
			\draw (-4.9,0.13) node {$-\infty$};
			\draw (-1,0.3) node {$-1$};
			\draw [domain=-5:-1.0,samples=100] plot(\x, {
				-ln(1-\x)/4+1/(2*(1-\x)) +ln(-1-\x)/4
			} );
		\end{tikzpicture}
		\caption{Function $P_{1}$}
	\end{minipage}
\begin{minipage}[t]{0.49\linewidth}
\centering
\begin{tikzpicture}[scale = 1]
	\draw (-2,0)--(2,0);
	\draw (0,2) -- (0,-2);
	\draw[dashed] (-1,-2) -- (-1,2);
	\draw[dashed] (1,-2) -- (1,2);
	\draw (-1.3,-0.2) node {$-1$};
	\draw (1.18,-0.2) node {$1$};
	\draw [domain=-0.99:0.7]
	plot(\x, {
		-ln(1-\x)/4+1/(2*(1-\x)) +ln(1+\x)/4
	} );
	samples = 100
\end{tikzpicture}
\caption{Function $P_{2}$}
	\end{minipage}
	\end{figure}
\begin{figure}[h]
\begin{center}
\begin{tikzpicture}[scale = 1]
	\draw (-0.6,0)--(3,0);
	\draw (0,-2.8) -- (0,0.6);
	\draw[dashed] (1,-2.8) -- (1,0.6);
	\draw (3,0.13) node {$+\infty$};
	\draw (1.14,-0.19) node {$1$};
	\draw [domain=1.15:2.8,samples=100] plot(\x, {
		-ln(-1+\x)/4+1/(2*(1-\x)) +ln(1+\x)/4
	} );
\end{tikzpicture}
		\caption{Function $P_{3}$}
	\end{center}
\end{figure}
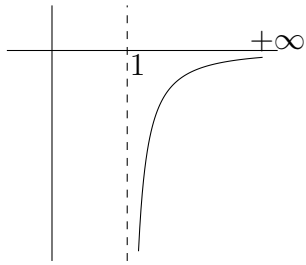

\newpage

Now, we are going to construct the inverse function of $P_1:(-\infty,-1)\rightarrow\mathbb{R}$. Recall $P_{1}\left(w(s)\right)=s-s_0\in(-\infty,0)$ so that $s<s_0$. Therefore, given $s_0\in\mathbb{R}$, we define
\begin{align*}
	t_1:(-\infty,s_0)\rightarrow(-\infty,-1),\, s\mapsto t_1(s)=P_1^{-1}(s-s_0).
\end{align*}
We can now define $f_2$,
\begin{align*}
	f_2:(-\infty,s_0)\rightarrow\mathbb{R},\, s\mapsto f_2(s)=\int_{s_0}^{s}t_1(v)dv+\widehat{f_2},\, \widehat{f_2}\in\mathbb{R}.
\end{align*}
We construct the inverse function of  $P_2:(-1,+1)\rightarrow\mathbb{R}$. Recall $P_2\left(w(s)\right)=s-s_0\in(-\infty,+\infty)$. Therefore, given $s_0\in\mathbb{R}$, we define
\begin{align*}
	t_2:\mathbb{R}\rightarrow(-1,1),\, s\mapsto t_2(s)=P_2^{-1}(s-s_0).
\end{align*}
We can define $f_3$,
\begin{align*}
	f_3:\mathbb{R}\rightarrow\mathbb{R},\, s\mapsto f_3(s)=\int_{s_0}^{s}t_2(r)dr+\widehat{f_3},\,\widehat{f_3}\in\mathbb{R}.
\end{align*}
Similarly,  given $s_0\in\mathbb{R}$, we define
\begin{align*}
	t_3:\left(-\infty,s_0\right)\rightarrow(1,+\infty),\, s\mapsto t_3(s)=P_3^{-1}(s-s_0).
\end{align*}
We can define $f_4$,
\begin{align*}
	f_4:\left(-\infty,s_0\right)\rightarrow\mathbb{R},\, s\mapsto f_4(s)=\int_{s_0}^{s}t_3(k)dk+\widehat{f_4},\,\widehat{f_4}\in\mathbb{R}.
\end{align*}

\subsubsection{Assume that $m>1$}
We have from \eqref{5}
\begin{equation}
\mathbf{Q(x)}:=\int\frac{dx}{(1-x^2)(1-mx)}=\frac{1}{2m-2}\ln\lvert1-x\lvert+\frac{1}{2m+2}\ln\lvert1+x\lvert-\frac{m}{m^2-1}\ln\lvert1-mx\lvert.
\end{equation}
For $m>1,$ it is clear that we need the following intervals
\begin{equation*}
	\mathbb R-\left\{-1, +1,1/m\right\}=(-\infty,-1)\cup\left(-1,1/m\right)\cup\left(1/m,+1\right)\cup(+1,+\infty),
\end{equation*} and we obtain different functions, by restricting function $\mathbf{Q}$ to the different intervals, namely, $Q_1=\mathbf{Q}\vert_{(-\infty,-1)}$, $Q_2=\mathbf{Q}\vert_{\left(-1,1/ m\right)}$, $Q_3=\mathbf{Q}\vert_{\left(1/m,+1\right)}$, $Q_4=\mathbf{Q}\vert_{(+1,+\infty)}$. We show some properties of them:

\begin{enumerate}
	\item $\lim_{\begin{subarray}{l}
			x\rightarrow-1\\
			x<-1
	\end{subarray}}Q_1(x)=-\infty$, $\lim_{x\rightarrow-\infty}Q_1(x)=-\frac{m}{m^2-1}\ln(m)<0$; decreasing, $\mathrm{Im}(Q_1)=\left(-\infty,\frac{m}{1-m^2}\ln(m)\right)$.
\item $\lim_{\begin{subarray}{l}
		x\rightarrow-1\\
		x>-1
\end{subarray}}Q_2(x)=-\infty$, $\lim_{\begin{subarray}{l}
		x\rightarrow\frac{1}{m}\\
		x<\frac{1}{m}
\end{subarray}}Q_2(x)=+\infty$; $Q_2(x)$ increasing, $\mathrm{Im}(Q_2)=\mathbb{R}$.
\item $\lim_{\begin{subarray}{l}
		x\rightarrow\frac{1}{m}\\
		x>\frac{1}{m}$, $
\end{subarray}}Q_3(x)=+\infty$, $\lim_{\begin{subarray}{l}
		x\rightarrow+1\\
		x<1
\end{subarray}}Q_3(x)=-\infty$; $Q_3(x)$ decreasing, $\mathrm{Im}(Q_3)=\mathbb{R}$.
\item $	\lim_{\begin{subarray}{l}
		x\rightarrow+1\\
		x>1
\end{subarray}}Q_4(x)=-\infty$, $\lim_{x\rightarrow+\infty}Q_4(x)=\frac{m}{1-m^2}\ln(m)<0$; $Q_4(x)$ increasing, $\mathrm{Im}(Q_4)$$=\left(-\infty,\frac{m}{1-m^2}\ln(m)\right)$.
\end{enumerate}

\begin{figure}[h!]
\begin{minipage}[t]{0.49\linewidth}
\centering
\begin{tikzpicture}[scale = 1]
	\draw (-5,0)--(0.5,0);
	\draw (0,0.5) -- (0,-3);
	\draw[dashed] (-1,0.5) -- (-1,-3);
	\draw[dashed] (-5,-0.42) -- (0.5,-0.42);
	\draw (-5,0.2) node {$-\infty$};
	\draw (-0.8,0.1) node {$-1$};
	\draw (0.96,-0.79) node {$\frac{m}{1-m^2}\ln(m)$};
	\draw [domain=-5:-1.01,samples=100] plot(\x, {
		ln(1-\x)/2+ln(-1-\x)/6-2*ln(1-2*\x)/3
	} );
	samples = 100
\end{tikzpicture}
\caption{Function $Q_{1}$}
\end{minipage}
\begin{minipage}[t]{0.49\linewidth}
	\centering
\begin{tikzpicture}[scale = 1]
	\draw (-2,0)--(2,0);
	\draw (0,2) -- (0,-2);
	\draw[dashed] (-1,-2) -- (-1,2);
	\draw[dashed] (0.5,-2) -- (0.5,2);
	\draw (-1.28,-0.2) node {$-1$};
	\draw (0.9,-0.2) node {$1/m$};
	\draw [domain=-0.99:0.49] plot(\x, {
		ln(1-\x)/2+ln(1+\x)/6-2*ln(1-2*\x)/3
	} );
	samples = 100
\end{tikzpicture}
\caption{Function $Q_{2}$}
\end{minipage}
\end{figure}
		
\begin{figure}[h!]
	\begin{minipage}[t]{0.49\linewidth}
		\centering
		\begin{tikzpicture}[scale = 1]
			\draw (-1,0)--(3,0);
			\draw (0,2.1) -- (0,-2);
			\draw[dashed] (1,-2) -- (1,2.1);
			\draw[dashed] (0.5,-2) -- (0.5,2.1);
			\draw (1.23,-0.17) node {$1$};
			\draw (0.3,-0.27) node {$\frac{1}{m}$};
			\draw [domain=0.51:0.99] plot(\x, {
				ln(1-\x)/2+ln(1+\x)/6-2*ln(-1+2*\x)/3
			} );
			samples = 100
		\end{tikzpicture}
		\caption{Function $Q_{3}$}
	\end{minipage}
	\begin{minipage}[t]{0.49\linewidth}
		\centering
		\begin{tikzpicture}[scale = 1]
			\draw (-0.99,0)--(3,0);
			\draw (0,1.2) -- (0,-3);
			\draw[dashed] (1,1.1) -- (1,-3);
			\draw[dashed] (-0.99,-0.45) -- (3,-0.45);
			\draw (3,0.1) node {$+\infty$};
			\draw (1.16,-0.19) node {$1$};
			\draw (-0.92,-0.7) node {$\frac{m}{1-m^2}\ln(m)$};
			\draw [domain=1.01:3] plot(\x, {
				ln(-1+\x)/2+ln(1+\x)/6-2*ln(-1+2*\x)/3
			} );
		\end{tikzpicture}
		\caption{Function $Q_{4}$}
	\end{minipage}
\end{figure}

\newpage
\begin{theorem}\label{Theorem 3.1}
Consider a space-like translator $\Gamma:\Omega\rightarrow\mathbb{H}^{n}\times_{-1}\mathbb{R}$, $\Gamma(p)=\left(p,u(p)\right)$ for all $p\in\Omega$, where $u\in C^{2}\left(\Omega\right)$. Then, $\Gamma$ is foliated by horospheres if, and only if, there exists $f\in C^{2}\left(I\right)$, $I\subseteq\mathbb{R}$ such that $u=f\circ\rho$ and $f$ is one of the following: 
	\begin{enumerate}
		\item $f_{1}:\mathbb{R}\rightarrow\mathbb{R},\, f_{1}(s)=\frac{s}{n-1}+f_{0},\, f_{0}\in\mathbb{R}$.
		\item $n=2,\, f_{3}:\mathbb{R}\rightarrow\mathbb{R},\, f_{3}(s)=\int {P_{2}}^{-1}(s)ds+f_0,\, f_{0}\in\mathbb{R}$.
		\item $n>2,\, f_{6}:\mathbb{R}\rightarrow\mathbb{R},\, f_{6}(s)=\int{Q_{2}}^{-1}(s)ds+f_0,\, f_{0}\in\mathbb{R}$.
		\item $n>2,\, f_{7}:\mathbb{R}\rightarrow\mathbb{R},\, f_{7}(s)=\int{Q_{3}}^{-1}(s)ds+f_0,\, f_{0}\in\mathbb{R}$.
	\end{enumerate}
\end{theorem}
\begin{theorem}\label{Theorem 3.2}
	Consider a graphical time-like translator $\Gamma:\Omega\rightarrow\mathbb{H}^{n}\times_{-1}\mathbb{R}$, $\Gamma(p)=\left(p,u(p)\right)$ for all $p\in\Omega$, where $u\in C^{2}\left(\Omega\right)$. Then, $\Gamma$ is foliated by horospheres if, and only if, there exists $f\in C^{2}\left(I\right)$, $I\subseteq\mathbb{R}$ such that $u=f\circ\rho$ and $f$ is one of the following: 
	\begin{enumerate}
		\item $n=2,\, f_{2}:\left(-\infty,s_0\right)\rightarrow\mathbb{R},\, f_{2}(s)=\int {P_{1}}^{-1}(s)ds+f_{0},\, f_{0}\in\mathbb{R}$.
		\item $n=2,\, f_{4}:\left(-\infty,s_0\right)\rightarrow\mathbb{R},\, f_{4}(s)=\int {P_{3}}^{-1}(s)ds+f_0,\, f_{0}\in\mathbb{R}$.
		\item $n>2,\, f_{5}:\left(-\infty,\frac{1-n}{n^2-2n}\ln(n-1)\right)\rightarrow\mathbb{R},\, f_{5}(s)=\int{Q_{1}}^{-1}(s)ds+f_0,\, f_{0}\in\mathbb{R}$.
		\item $n>2,\, f_{8}:\left(-\infty,\frac{n-1}{2n-n^2}\ln(n-1)\right)\rightarrow\mathbb{R},\, f_{8}(s)=\int{Q_{4}}^{-1}(s)ds+f_0,\, f_{0}\in\mathbb{R}$.
	\end{enumerate}
\end{theorem}
\section{Rotationally Invariant Translators}\label{section:space-like}
We introduce the \textit{Lorentzian Space}
\begin{equation*}
	\mathbb{L}^{n+1}=\mathbb{R}^{n+1},\, {\langle\left(x_{1},\ldots,x_{n+1}\right), \left(y_{1},\ldots,y_{n+1}\right)\rangle}_{L}=x_{1}y_{1}+\ldots+x_{n}y_{n}-x_{n+1}y_{n+1},\, n\geq2,
\end{equation*} We consider the model of the hyperbolic space
\begin{equation*}
	\mathbb{H}^{n}=\left\{p\in\mathbb{L}^{n+1}/\,{\langle p,p \rangle}_{L}=-1,\, p_{n+1}>0\right\}.
\end{equation*}
Simple computations provide $T_{p}\mathbb{H}^{n}=\left\{V\in\mathbb{R}^{n+1}/{\langle V,p\rangle}_{L}=0\right\}$.
Now, let us define
$\chi:\mathbb{H}^{n}\rightarrow\mathbb{R}^{n+1}$, $\chi(p)=p$ \textit{the position vector}, which is a unit normal time-like vector field to $\mathbb{H}^{n}$ on $\mathbb{L}^{n+1}$.
Since $A_{\chi}V=+V,$ equivalently $A_{\chi}=+{id}_{T_{p}\mathbb{H}^{n}},$ it is clear that $\mathbb{H}^{n}\rightarrow\mathbb{L}^{n+1}$ is \textit{totally umbilical}.
We define
$$\rho:\mathbb{H}^{n+1}\rightarrow[1,+\infty),\,\left(p_{1},\ldots,p_{n}\right)\mapsto p_{n+1}.$$
We have $\rho^{-1}\left\{s\right\}=\left\{p\in\mathbb{H}^{n}\vert p_{n+1}=s,\ s\geq1\right\},$ so that this set is a $(n-1)$ $\mathrm{sphere}$ of radius $\sqrt{s^{2}-1}$. Simple computations show
\begin{gather*}
	\nabla\rho=-\partial_{n+1}-{\langle\partial_{n+1},\chi\rangle}_{_L}\chi, \\
	{\langle\nabla\rho,\nabla\rho\rangle}_{L}=\langle-\partial_{n+1}-\langle\partial_{n+1},\chi\rangle\chi,-\partial_{n+1}-\langle\partial_{n+1},\chi\rangle\chi\rangle=\langle\partial_{n+1},\chi\rangle^{2}-1 \neq+1.
\end{gather*}
We reparametrize the projection, obtaining
\[ \tau:\mathbb{H}^n\to(0,+\infty), \quad \tau=\ln\left(\sqrt{\rho^{2}-1}+\rho\right).
\]
Clearly, $\rho=\cosh(\tau)$, and now $\nabla\tau=\frac{1}{\sqrt{\rho^{2}-1}}\nabla\rho$, is a unit normal to each level set. We compute the mean curvature of each level set, namely,
\begin{equation*}
	h=\mathrm{div}(\nabla\tau)=\mathrm{div}\left(\frac{\nabla\rho}{\sqrt{\rho^{2}-1}}\right)=\left(n-1\right)\coth(\tau).
\end{equation*} Assume that $u=f\circ\rho\,$ for some $f:I\rightarrow\mathbb{R}$. This means that $\Gamma_{u}$ is foliated by spheres. According to \textit{Theorem 3.5} of paper \cite{LO}, we have:
\begin{proposition}
	 A function $u=f\circ\tau$  provides a graphical, rotationally invariant, translator $\Gamma_u$ if, and only if, the function $f:I\subset (0,+\infty)\to\R$ is a solution to the following ODE:
	\begin{equation}\label{2.3}
		f''(s)=\left(1-(f'(s))^{2}\right)\left(1-(n-1)\coth(s)f'(s)\right).
	\end{equation}
\end{proposition} Our next target is to study the solutions to \eqref{2.3}. For this aim, we take $w=f'$ in \eqref{2.3} and deal with the following ODE
\begin{equation}\label{2.4}
	w'(s)=\left(1-w^{2}(s)\right)\left(1-(n-1)\coth(s)w(s)\right).
\end{equation}
\begin{remark} Given a solution $w$ to this ODE, each primitive $f=\int w$ will provide a graphical rotationally invariant translator in the following way. Define $u:=f\circ\tau$, being $\Gamma_u$ its graph.
\end{remark}

Firstly, we have the trivial solutions
\[ w_{+1}, w_{-1}:(0,+\infty)\to\R, \quad w_{+1}(s)=+1, \quad w_{-1}(s)=-1.
\]
The primitives $\mathbf{f}_{\pm 1}=\int w_{\pm 1}$ will provide degenerate hypersurfaces, but we will need them in the computations of rest of this section.

From now, we will use the book \cite{Wiggins}. Indeed, we define the following dynamical system:
\begin{equation}
	X:\mathbb{R}^{2}\rightarrow\mathbb{R}^{2},\,X(s,z)=\left(\sinh(s),\left(1-z^{2}\right)\left(\sinh(s)-(n-1)\cosh(s)z\right)\right).
\end{equation}
The zeros of $X$ are
\[p_{0}=(0,0),\, p_{1}=(0,1),\, p_{-1}=(0,-1).\]
We need to compute
\begin{align*}
	&\frac{\partial{X}}{\partial{s}}=\left(\cosh(s),(1-z^{2})\left(\cosh(s)-(n-1)\sinh(s)z\right)\right),\\
	&\frac{\partial{X}}{\partial{z}}=\left(0,-2z\left(\sinh(s)-(n-1)\cosh(s)z\right)+\left(1-z^{2}\right)\left((1-n)\coth(s)\right)\right).
\end{align*} We denote the differential of $X$ at $p$ by $DX(p)$. We classify the points $p_{0},\, p_{1}$ and $p_{-1}$ according to the eigenvalues of $DX(p)$:
\begin{equation*}
	DX(p_{0})=\begin{pmatrix}
		1 & 0 \\
		1 & 1-n
	\end{pmatrix},\ DX(p_{1})=\begin{pmatrix}
		1 & 0 \\
		0 & 2(n-1)
	\end{pmatrix},\ DX(p_{-1})=\begin{pmatrix}
		1 & 0 \\
		0 & 2(n-1)
	\end{pmatrix},\ n\geq2 .
\end{equation*}Clearly, $p_{1}$ and $p_{-1}$ are \textit{sources} because both eigenvalues are positive. But $p_{0}$ is a \textit{saddle point}. We introduce $\xi=\left\{(s,z)\in\mathbb{R}^{2}\vert s\geq0\right\}$,
$\xi_{+}=\left\{(s,z)\in\mathbb{R}^{2}\vert z>1\right\}$, $\xi_{-}=\left\{(s,z)\in\mathbb{R}^{2}\vert\right.$ $\left.z<-1\right\}$, and $\xi_{0}=\left\{(s,z)\in\mathbb{R}^{2}\vert\right.$ $\left.-1<z<1\right\}$.

\begin{lemma}\label{lemma 4.1}
	For each solution $w$ to \eqref{2.4}, there exist an integral curve of $X$.
\end{lemma}
\begin{proof}
	We take $w:I\rightarrow\mathbb{R}$ a solution to \eqref{2.4}. We consider a solution to $s'(r)=\sinh\left(s(r)\right)$, say $s(r)$. Take $z(r)=w\left(s(r)\right)$ so that $\gamma(r)=\left(s(r),z(r)\right)=\left(s(r),w\left(s(r)\right)\right)$ is an integral curve of $X$.
\end{proof}
\begin{lemma}\label{lemma 4.2}
	For each integral curve $\gamma:K_{o}\rightarrow\mathbb{R}^{2}$ such that, $\gamma(r)=\left(s(r),z(r)\right)$ with $s$ an bijective map, then $w:= z\circ s^{-1}$ is a solution to \eqref{2.4}.
\end{lemma}
\begin{proof}
	Since $\gamma$ is an integral curve of $X$ $\left(\text{that is,}\, X\left(\gamma(r)\right)=\gamma'(r)\right)$, then
	$$s'(r)=\sinh\left(s(r)\right),\quad z'(r)=\left(1-z(r)^{2}\right)\left(\sinh(s(r))-(n-1)\cosh(s(r))z(r)\right).$$
	As $s$ is bijective, then $s'\left(s^{-1}(y)\right)=\sinh(y)$ for any $y\in K_{o}$. Therefore,
	\begin{equation*}
		w'(y)=\frac{z'\left(s^{-1}(y)\right)}{s'\left(s^{-1}(y)\right)}=\left(1-w(y)^{2}\right)\left(1-(n-1)\coth(y)w(y)\right).
	\end{equation*} This completes the proof.
\end{proof}
\begin{lemma}\label{lemma 4.3}
	For each solution to the ODE $z'(r)=\left(1-n\right)\left(1-z(r)^{2}\right)z(r)$, the curve $\beta(r)=\left(0,z(r)\right)$ is an integral curve of $X$.
\end{lemma}
\begin{proof} We take the curve
	$\beta(r)=\left(0,z(r)\right)$. Then, we show that $\beta$ is an integral curve of $X$. Indeed,
	$\beta'(r)=\left(0,z'(r)\right)=X\left(0,z(r)\right)=\left(\sinh(0),(1-z(r)^{2})(\sinh(0)-(n-1)\cosh(0)\right.$ $\left.z(r))\right)$ $=\left(0,(1-z(r)^{2})(1-n)z(r)\right)$, that is, $z'(r)=(1-n)(1-z(r)^{2})z(r)$.
\end{proof}
\subsection{The Space-like Case}

The space-like translators appear when $\vert f'\vert<1$, that is, when $\vert w\vert <1$.

\begin{lemma}\label{lemma 4.4}
	Each solution w to \eqref{2.4} with initial condition $(s_{0},z_{0})\in\xi_{0}$, namely,   $w(s_{0})=z_{0}$, can be extended (as solution) to $w:(0,+\infty)\rightarrow\mathbb{R}$. Also, $\lim_{s\rightarrow0}w(s)\in\left\{-1,0,+1\right\}$.
\end{lemma}
\begin{proof}
	We take a local solution $w$ to \eqref{2.4} such that $-1<w(s)<+1$, which is bounded by the constant solutions $w_{\pm}(s)=\pm1$. Therefore, $w$ can be extended to $w:\left(0,+\infty\right)\rightarrow\mathbb{R}$. By Lemma \ref{lemma 4.1}, we construct an integral curve $\gamma_{w}:\left(0,+\infty\right)\rightarrow\mathbb{R}$ of $X$ from $w$. By Lemma \ref{lemma 4.3}, $\gamma_{w}$ and $\beta$ can only coincide at some point $p$ such that $X(p)=0$. But this only holds when $p\in\left\{p_{0},p_{+1},p_{-1}\right\}$, namely, when $\lim_{s\rightarrow0}w(s)\in\left\{-1,0,+1\right\}$.
\end{proof}
We define of function $\vartheta:\xi\rightarrow\mathbb{R}$ given by
\begin{equation}\label{theta}
	\vartheta(s,z)=\left(1-z^{2}\right)\left(\sinh(s)-(n-1)\cosh(s)z\right).
\end{equation} We need the set $\alpha=\left\{(s,z)\in\xi \mid\vartheta(s,z)=0\right\}=\left\{(s,z)\in\xi \mid z=\tanh(s)/(m-1)\right\}$.

\begin{proposition}{\label{proposition4.2}}
	Let $w:(0,+\infty)\to (-1,1)$ be a solution to \eqref{2.4} such that $\lvert w(s_{0})\rvert<1$  for some $s_{0}\in\left(0,+\infty\right)$.  Then, it is one of the following:
	\begin{enumerate}
		\item There exists a unique $\hat{w}:[0,+\infty)\rightarrow(-1,+1)$ solution to \eqref{2.4} such that $\hat{w}(0)=0$. For all $s\in(0,+\infty)$, it holds $\hat{w}(s)<\frac{\tanh(s)}{n-1}$.
		\item If  $w(s_0)<\tanh(s_0)/(n-1)$, then
		$w(s)<\tanh(s)/(n-1)$ for any $s\ge s_0$. In addition, if $w(s_0)<\hat{w}(s_0)$, then
		$\lim_{s\rightarrow0}w(s)=-1$.
		\item If $w(s_{0})>\hat{w}(s_{0})$, then $\lim_{s\rightarrow0}w(s)=+1$.
	\end{enumerate}
	In all cases, it holds $\lim_{s\rightarrow+\infty}w(s)=\frac{1}{n-1}$.
\end{proposition}
\begin{proof}
	
	\textit{1.} We want to solve the boundary problem:
	\[	w'(s_{0})=\left(1-w(s)^{2}\right)\left(1-(n-1)\cosh(s)w(s)\right),\quad w(0)=0.
	\]
	If there is a solution, the associated integral curve of $X$ will start at $(0,0)$.
	The eigenvalues of $DX(0,0)$ are $\lambda_{1}=1,\lambda_{2}=1-n$ and then their eigenvectors are $V_{1}=(n,1)$ and $V_{2}=(0,1)$ respectively. By \textit{Theorem 3.2.1.} of the book \cite{Wiggins}, there exists an integral curve $\gamma$ of $X$ such that $\gamma(0)=(0,0)$ and $\gamma'(0)=V_{1}$ (Integral submanifolds). Then, by Lemma \ref{lemma 4.2}, we obtain the needed solution to the boundary problem.
	
	In addition, assume by contradiction that there exist $s_{0}\in(0,+\infty)$ such that $\hat{w}(s_{0})>\frac{\tanh(s)}{n-1}>0$. We have $\hat{w}'(s_{0})<0$ and $\hat{w}(0)=0$. By continuity, there exist $s_{1}\in(0,s_{0})$ such that $\hat{w}(s_{1})=\frac{\tanh(s_{1})}{n-1}$. There exist $s_{2}\in(s_{1},s_{0})$ such that $\hat{w}'(s_{1})>0$ and $\hat{w}'(s_{2})>\frac{\tanh(s_{2})}{n-1}$. But then,
	\begin{align*}
		&0<\hat{w}'(s_{2})=\left(1-{\hat{w}'(s_{2})}^{2}\right)\left(1-(n-1)\coth(s_{0})\hat{w}(s_{2})\right)<0,
	\end{align*}which is a contradiction.
	
	\textit{2.} We recall function $\vartheta$, \eqref{theta}. Assume that $(s_0,w(s_0))$ satisfies $\vartheta(s_0,w(s_0))=0$. Then, $w'(s_0)=(1-w(s_0)^2)\vartheta(s_0,w(s_0))/\tanh(s_0)=0$. Moreover, there exists  $\delta>0$ small enough such that if $s\in(s_0-\delta,s_0)$, then $w'(s)=(1-w(s_0)^2)\vartheta(s_0,w(s_0))/\tanh(s_0)<0$; and if $s\in(s_0,s_0+\delta)$, then $w'(s)>0$. That is to say, if the graph of $w$ touches the set $\alpha$, then $w$ has a local minimum. Therefore, if the graph of $w$ is below $\alpha$ at some point, then it remains below $\alpha$ for any other further point.
	
	In addition, if $w(s_0)<\hat{w}(s_0)$, by Lemma \ref{lemma 4.4}, $\lim_{s\rightarrow0}w(s)\in\left\{-1,0,+1\right\}$. Since $z_{0}<\hat{w}(s_{0})$, by uniqueness of solutions to ODE, we have $w(s)<\hat{w}(s)$ for all $s\in\left(0,+\infty\right)$. Also, $\hat{w}(0)=0>\lim_{s\rightarrow0}w(s)$ and therefore, $\lim_{s\rightarrow0}w(s)=-1$.
	
	\textit{3.} If $w(s_{0})>\hat{w}(s_{0})$, then $w'(s_0)=(1-w(s_0)^2)\vartheta(s_0,w(s_0))/\tanh(s_0)<0$. Therefore, the graph of $w$ remains above  the graph of $\hat{w}$ in $(0,s_0)$. By Lemma \ref{4.2}, the only possible limit is $\lim_{s\to 0}w(s)=+1$.
	
	Finally, when $s$ becomes big, if the graph of $w$ is above $\alpha$, then $w$ is strictly decreasing. However, if the graph of $w$ is below $\alpha$, $w$ is strictly increasing. Then, $\lim_{s\to +\infty}w(s)=\lim_{s\to+\infty} \tanh(s)/(n-1)=1/(n-1)$.
\end{proof}
\begin{definition}
	We define the function $\hat{f}(s)=\int\hat{w}(s)ds+f_{0}$. The associated rotationally invariant \textit{translator} is called \textit{the bowl}.
\end{definition}
By Proposition \ref{proposition4.2}, we immediately have the following result.
\begin{theorem}\label{theorem 4.1}
	Given a rotationally symmetric space-like translator $\Gamma$ in $\mathbb{H}^{n}\times_{-1}\mathbb{R}$, there exist $f:(0,+\infty)\rightarrow\mathbb{R}$ such that $f$ is a solution to the ODE \eqref{2.4}, and  $\Gamma(p)=\left(p,f(\rho(p))\right)$, for all $p\in\mathbb{H}^{n}$. There are 3 types of functions $f=\int w$, namely the bowl, those with $\lim_{s\to 0}w(s)=1$ and those with $\lim_{s\to 0}w(s)=-1$.
\end{theorem}
\begin{remark}
In items 2 and 3 Proposition \ref{proposition4.2}, $\lim_{s\rightarrow 0}w(s)=\pm1$. This means that when approaching to the axis of rotation, the translator will hit the axis with an angle of $\pi/4$. Therefore, there is a conic singularity.
\end{remark}
\subsection{The Time-like Case}
\begin{lemma}\label{lemma 4.10}
	Given $\left(s_{0},z_{0}\right)\in\mathbb{R}^{2}$ such that $s_{0}>0$, $z_{0}>1$, the associated solution can be extenden to $w:[0,s_{0}+\varepsilon)\rightarrow[1,+\infty)$ with $\lim_{s\rightarrow0}w(s)=+1$.
\end{lemma}
\begin{proof}
	The constant solution $w_{+}(s)=+1$ is a bound from below. Since, $\vartheta(s,z)>0$, for all $\left(s,z\right)\in\xi_{+}$. We know $w$ is increasing. Then, we can extend to $[0,s_{0}+\varepsilon)$ similarly to item 2 of Proposition \ref{proposition4.2}, $\lim_{s\rightarrow0}w(s)=+1$.
\end{proof}
\begin{lemma}\label{lemma 4.6}
	Given $\left(s_{0},z_{0}\right)\in\mathbb{R}^{2}$ such that $s_{0}<0$, $z_{0}<1$, the associated solution can be extended to $w:[0,s_{0}+\varepsilon)\rightarrow(-\infty,-1]$ with $\lim_{s\rightarrow0}w(s)=-1$. Since, $\vartheta(s,z)<0$, forall $\left(s,z\right)\in\xi_{-}$.
\end{lemma}
\begin{proof}
	The proof is very similar to Lemma \ref{lemma 4.10}.
\end{proof}
Take $g$ a function which is the inverse of $f$. We compute $1=g'\left(f(s)\right)f'(s)$ and therefore $0=g''\left(f(s)\right)f'(s)^{2}+g'\left(f(s)\right)f''(s)$. Then, $$g''\left(f(s)\right)=-\frac{\left(1-f'(s)^{2}\right)\left(1-\coth(s)f'(s)\right)}{f'(s)^{3}}.$$ We take $f\circ g(t)=t$, which implies $0=f''\left(g(t)\right)g'(t)^{2}+f'\left(g(t)\right)g''(t)$. We obtain the following ODE
\begin{equation}\label{16}
	g''(t)=\left(g'(t)-1\right)\left(\coth\left(g(t)\right)-g'(t)\right).
\end{equation}
\begin{lemma}\label{lemma 4.12}
	Given a solution $g:\left(t_{0}-\varepsilon,t_{0}+\varepsilon\right)\rightarrow\mathbb{R}$ to ODE \eqref{16}. Then, $g''(t_{0})=-\coth(g_{0})$  $<0$, there are two functions $f_{\pm}:(g_{0}-\delta,g_{o}]\rightarrow\mathbb{R}$ solutions to the \eqref{2.3}, which are inverse functions of $g$ with $\lim_{s\rightarrow g_{0}}f'_{\pm}(s)=\lim_{s\rightarrow g_{0}}\frac{1}{g'\left(f_{\pm}(s)\right)}=\pm\infty$. We construct two graphical translator from $f_{\pm}$ rotationally invariant and together they make a smooth hypersurface.
\end{lemma}
\begin{proof}
	Choose $t_{0}\in\mathbb{R}$. We consider the following IVP:
	\begin{equation*}
		g''(t)=\left(\varepsilon'+\tilde{\varepsilon}g'(t)^{2}\right)\left(\coth\left(g(t)-g'(t)\right)\right),\quad g'(t_{0})=0,\,g(t_{0})=s_{0}\in I.
	\end{equation*} As usual, there exists a smooth solution $\alpha:\left(t_{0}-\varepsilon,t_{0}+\varepsilon\right)\rightarrow\mathbb{R}$. Note that $t_{0}$ is a critical point of $\alpha$ and $\alpha''(t_{0})=\varepsilon'\coth(s_{0})$.
	Case $\coth(s_{0})\neq0:$ Then, $t_{0}$ is an extremum of $\alpha$. The restrictions $g_{+}=g\vert\left(t_{0},t_{0}+\varepsilon\right)$ and $g_{-}=g\vert\left(t_{0}-\varepsilon,t_{0}\right)$ will be injective, by reducing $\varepsilon$ if necessary. Construct their inverse functions $f_{+}$ and $f_{-}$ satisfy \eqref{2.3}. To do so, we put $f_{+}\left(g(t)\right)=t$, and therefore
	\begin{equation*}
		\begin{aligned}
			1&=f'_{+}\left(g(t)\right)g'(t), \quad 0=f''_{+}\left(g(t)\right)g'(t)^{2}+f'_{+}\left(g(t)\right)g''(t),\\
			f''_{+}\left(g(t)\right)g'(t)^{2}&=-f'_{+}\left(g(t)\right)\left(\varepsilon'+\tilde{\varepsilon}\alpha'(t)^{2}\right)\left(\coth\left(g(t)-g'(t)\right)\right).
		\end{aligned}
	\end{equation*} Next, we change $s=g(t)$, and then $g'(t)=1/f'_{+}(s)$, so that
	\begin{equation*}
		\frac{f''_{+}(s)}{f'_{+}(s)^{2}}=\frac{1}{f'_{+}(s)^{2}}\left(\tilde{\varepsilon}+\varepsilon'f'_{+}(s)^{2}\right)\left(1-\coth(s)f_{+}(s)\right).
	\end{equation*} A similar computation holds for $f_{-}$. The union of the corresponding graphical translators and their common boundary provide a smooth translator, because $g$ is a smooth map and $f_{+}$, $f_{-}$ are tools to reparametrize its graph.
\end{proof}
\begin{definition}
	The spindle is the rotationally invariant hypersurface obtained in Lemma \ref{lemma 4.12}.
\end{definition}
\begin{remark}
According to Lemmas \ref{lemma 4.10} and \ref{lemma 4.6}, $\lim_{s\rightarrow 0}w(s)=\pm1$. Again, we have a conic singularity  of angle $\pi/4$.
\end{remark}
\begin{theorem}\label{4.2}
	Any rotationally invariant, time-like translator in $\mathbb{H}^{n}\times_{-1}\mathbb{R}$ is an open subset of a spindle.
\end{theorem}
\begin{proof}
	Take $w:\left(s_{0}-\varepsilon,s_{0}+\varepsilon\right)\rightarrow\left(1,+\infty\right)$ a solution to \eqref{2.4}. By Lemma \ref{lemma 4.10}, we can extend $w:\left[0,s_{0}+\varepsilon\right)\rightarrow\left(1,+\infty\right)$, and $\lim_{s\rightarrow0}w(s)=1$. We will show that there exist $s_{1}>s_{0}+\varepsilon$ such that $\lim_{s\rightarrow s_{1}}w(s)=+\infty$ (finite-time blow-up). We have $$w'(s)=(1-w(s)^2)(1-(n-1)\coth(s)w(s))=\sinh(s)\vartheta(s,w(s)).$$ We can get $F(s,z)=\left(1-z^{2}\right)\left(1-(n-1)\coth(s)z\right)$. With this, we have $w'(s)=F\left(s,\right.$ $\left.w(s)\right)$. We define $\Upsilon=\left(0,+\infty\right)\times\left[1,+\infty\right)$, and  $F,G:\Upsilon\rightarrow\mathbb{R}$, we consider $G(s,z)=\left(z^{2}-1\right)z$. Clearly, $F(s,z)\geq G(s,z)$, for any $(s,z)\in\Upsilon$. Next, the solution to $z'(s)=G\left(s,z(s)\right)$ with $z\left(s_{0}\right)=z_{0}$ is $z:\left(0,s_{0}+A\right)\rightarrow\mathbb{R}$, $$z(s)=\frac{1}{\sqrt{1-e^{2s-2(s_{0}+A)}}}\,,\quad A=-\frac{1}{2}\ln\left(1-\frac{1}{{z_{0}}^{2}}\right),\, A\in\mathbb{R}.$$
	Note that $\lim_{s\rightarrow s_{0}+A}z(s)=+\infty$. These previous computations imply that there exist $s_{1}\in\left(s_{0},s_{0}+A\right)$ such that $\lim_{s\rightarrow s_{1}}w(s)=+\infty$. We start with $f$, and then we take $w=f'>1$. Every $w$ has a finite-time blow-up, so we use Lemma \ref{lemma 4.12}. Therefore every $f$ is strictly increasing and $g=f'$ so at some point, $\lim_{z\rightarrow z_{0}}g'(z)=0$.  The same reasoning works when $w=f'<-1$.
\end{proof}
\section{Isometries and a Quasilinear Elliptic PDE}\label{section:time-like}
\begin{lemma}\label{lemma 5.1}
When $\varepsilon=-1$, PDE \eqref{PDE 2} behaves as a quasilinear elliptic operator. Moreover, it is locally uniformly elliptic.
\end{lemma}
\begin{proof}
We use the model of half-plane :
\begin{equation*}
\begin{aligned}
\mathbb{H}^{n}=\left\{\left(x_1,\ldots,x_n\right)\in \mathbb{R}^{n}\vert x_1>0 \right\},\, n\geq2.\\
B_{1}=\left(\partial_{1},\ldots,\partial_{n}\right),\, B_{2}=\left(e_{i}=x_{1}\partial_{i} : i=1,\ldots,n\right).
\end{aligned}
\end{equation*} We know the hyperbolic metric is $g=g_{\mathbb{H}^{n}}=\frac{1}{{x_1}^{2}}\langle,\rangle\,;\, g\left(e_{i},e_{j}\right)={x_{1}}^{2}\, g\left(\partial_{i},\partial_{j}\right)=\langle\partial_{i},\partial_{j}\rangle=\delta_{ij}$, where as usual, $\delta_{ij}$ denotes the Kronecker's delta. We take $u_{i}=\partial_{i}u$, and then $$\nabla u=x_{1}^{2}\sum_{i=1}^{n}u_{i}\partial_{i}.$$ Take $\lvert\nabla u \rvert ^{2}=\sum_{i=1}^{n}{x_{1}}^{2}{u_{i}}^{2};\, \nabla_{\partial_{i}}\partial_{j}=\sum_{k=1}^{n}\Gamma_{ij}^{k}\partial_{k}$. We assume $\lvert\nabla u \rvert ^{2}<1$, then $1>{x_{1}}^{2}\sum_{i=1}^{n}$ ${u_{i}}^{2}$ . We make some computations :
\begin{equation*}
	\frac{1}{W}=\mathrm{div}\left(\frac{\nabla u}{W}\right)=\frac{-(\nabla u)(W)}{W^{2}}+\frac{1}{W}\mathrm{div}(\nabla u),
\end{equation*}
\begin{equation}\label{eq 15}
	1=\mathrm{div}(\nabla u)-\frac{(\nabla u)(W)}{W}.
\end{equation}
\begin{equation*}
\begin{aligned}
\mathrm{div}(\nabla u)&=\sum_{i=1}^{n}g\left(\nabla_{e_{i}}\nabla u,e_{i}\right)={x_{1}}^{2}\sum_{i=1}^{n}u_{ii}+2u_{1}x_{1}+\sum_{i,j=1}^{n}{x_{1}}^{2}\Gamma_{ij}^{i}u_{j}.\\
W&=\left(1-\lvert\nabla u \rvert ^{2}\right)^{1/2}=\left(1-\sum_{i=1}^{n}{x_{1}}^{2}{u_{i}}^{2}\right)^{1/2},
\end{aligned}
\end{equation*}
\begin{equation*}
(\nabla u)(W)=\left({x_{1}}^{2}\sum_{i=1}^{n}u_{i}\partial_{i}\right)\left(1-\sum_{j=1}^{n}{x_{1}}^{2}{u_{j}}^{2}\right)^{1/2}=-\frac{{x_{1}}^{3}}{W}u_{1}\sum_{j=1}^{n}{u_{j}}^{2}-\frac{{x_{1}}^{4}}{W}\sum_{i,j=1}^{n}u_{i}u_{j}u_{ij}.
\end{equation*}  We insert all these computations in \eqref{eq 15}.
\begin{equation*}
0=\sum_{i,j=1}^{n}\left(\delta_{ij}+\frac{{x_{1}}^{2}u_{i}u_{j}}{W^{2}}\right)u_{ij}+\frac{2u_{1}}{x_{1}}+\frac{x_{1}u_{1}}{W^{2}}\sum_{i=1}^{n}{u_{i}}^{2}+\sum_{i,j=1}^{n}u_{i}\Gamma_{ji}^{j}-\frac{1}{{x_{1}}^{2}}.
\end{equation*}
We multiply by $W^{2}x_{1}^{2}$. Later, we will that this is not a problem.
\begin{align*}
0&=\sum_{i,j=1}^{n}\left({x_{1}}^{2}\left(1-\sum_{k=1}^{n}{x_{1}}^{2}{u_{k}}^{2}\right)\delta_{ij}+{x_{1}}^{4}u_{i}u_{j}\right)u_{ij}+2u_{1}x_{1}\left(1-\sum_{k=1}^{n}{x_{1}}^{2}{u_{k}}^{2}\right)\\
&+{x_{1}}^{3}u_{1}\sum_{i=1}^{n}{u_{i}}^{2}+\sum_{i,j=1}^{n}{x_{1}}^{2}\left(1-\sum_{k=1}^{n}{x_{1}}^{2}{u_{k}}^{2}\right)u_{i}\Gamma_{ji}^{i}(x)-1+\sum_{k=1}^{n}{x_{1}}^{2}{u_{k}}^{2}.
\end{align*}
We define
\begin{equation*}
\begin{aligned}
&a^{ij}, b:\mathbb{H}^{n}\times\mathbb{R}^{n}\rightarrow\mathbb{R},\quad a^{ij}(x,p)={x_{1}}^{2}\left(1-\sum_{k=1}^{n}{x_{1}}^{2}{p_{k}}^{2}\right)\delta_{ij}+{x_{1}}^{4}p_{i}p_{j},\,\text{and}\\
b(x,p)&=2u_{1}x_{1}\left(1-\sum_{k=1}^{n}{x_{1}}^{2}{u_{k}}^{2}\right)+{x_{1}}^{3}u_{1}\sum_{i=1}^{n}{u_{i}}^{2}+\sum_{i,j=1}^{n}{x_{1}}^{2}\left(1-\sum_{k=1}^{n}{x_{1}}^{2}{u_{k}}^{2}\right)u_{i}\Gamma_{ji}^{i}(x)-1\\ &+\sum_{k=1}^{n}{x_{1}}^{2}{u_{k}}^{2}.
\end{aligned}
\end{equation*}
Let us check that the matrix $A=\left(a^{ij}\right)$ is positive-define.We rewrite it as
$$A={x_{1}}^{2}\left(1-\sum_{k=1}^{n}{x_{1}}^{2}{p_{k}}^{2}\right)I_{n}+{x_{1}}^{4}B,\, B=p^{t}p.$$
By taking $p,q\in\mathbb{R}^{n}$, we see $B=p^{t}p$, $\langle p,q\rangle =pq^{t}$,  $Bp^{t}=p^{t}pp^{t}=\lvert p\rvert^{2}p^{t}$, and therefore $$Ap^{t}=\left({x_{1}}^{2}\left(1-\sum_{k=1}^{n}{x_{1}}^{2}{p_{k}}^{2}\right)+{x_{1}}^{4}\lvert p\rvert^{2}\right)p^{t},$$ where $I_{n}$ is the identity matrix. Clearly, $Bp^{t}=\lvert p\rvert^{2}p^{t}$, and given $q\perp p$ such that $pq^{t}=0$, then, $Aq^{t}={x_{1}}^{2}\left(1-\sum_{k=1}^{n}{x_{1}}^{2}{p_{k}}^{2}\right)q^{t}$. Therefore, the eigenvalues of $A$ are:
$$\lambda_{1}={x_{1}}^{2}\left(1-\sum_{k=1}^{n}{x_{1}}^{2}{p_{k}}^{2}\right){x_{1}}^{4}\lvert p\rvert^{2},\quad \lambda_{2}(p)={x_{1}}^{2}\left(1-\sum_{k=1}^{n}{x_{1}}^{2}{p_{k}}^{2}\right).$$
We define
\begin{equation*}
\Lambda=\left\{(x,p)\in\Omega\times\mathbb{R}^{n}\vert\,{x_{1}}^{2}\lvert p\rvert^{2}<1,x_{1}>0\right\},\,\lambda_{1},\lambda_{2}:\Lambda\rightarrow\mathbb{R}.
\end{equation*}In $\Lambda$, $\lambda_{2}>0$ and $\lambda_{1}>0$, so that the operator $Q$ is elliptic. Also, given $\mathbb{U}\subset\Lambda$ such that $\mathbb{U}$ is open and $\bar{\mathbb{U}}$ compact, then in $\mathbb{U}$
$$\frac{\lambda_{1}}{\lambda_{2}}=1+\frac{{x_{1}}^{2}\lvert p\rvert^{2}}{1-\sum_{k=1}^{n}{x_{1}}^{2}{p_{k}}^{2}}\geq1,$$ is bounded from above and from below. That is to say, $Q$ is locally uniformly bounded.
\end{proof}
We recover \textit{Theorem 10.2} of paper \cite{libroEDP};
\begin{theorem}\label{Theorem 5.1}
Let $\Omega$ be a bounded open domain in $\mathbb{R}^{n}$. Let $u,v\in C^{0}\left(\bar{\Omega}\right)\cap C^{2}\left(\Omega\right)$. Let $Q$ be a quasilinear operator such that :
\begin{enumerate}
	\item $Q$ is locally uniformly elliptic with respect to either $u$ or $v$;
	\item $a^{ij}$ do not depend on $z$;
	\item $b$ is non-increasing in $z$ for each $(x,p)\in\Omega\times\mathbb{R}^{n}$;
	\item $a^{ij}$, $b$ are continuously differentiable with respect to the p-variables in $\Omega\times\mathbb{R}\times\mathbb{R}^{n}$.
\end{enumerate} Assume that $Qu=Qv$ in $\Omega$ and $u=v$ on $\partial\Omega$. Then $u\equiv v$ in $\Omega$.
\end{theorem}
From Lemma \ref{Theorem 5.1} and Theorem \ref{Theorem 5.1}, we immediately obtain the following uniqueness result.
\begin{theorem}\label{uniq}
Let $\Omega$ be a bounded open domain in $\mathbb{H}^{n}$. Let $u,v\in C^{0}\left(\bar{\Omega}\right)\cap C^{2}\left(\Omega\right)$ such that $\Gamma_{u}$, $\Gamma_{v}$ are space-like translators, and $u=v$ on $\partial\Omega$. 
Then $u\equiv v$ in $\Omega$.
\end{theorem}
\begin{proof}
By Lemma \ref{lemma 5.1}, Our operator $Q$ is quasilinear elliptic and locally uniformly bounded. In addition, $Q$ does not depend on $z$. We know $a^{ij},\, b\in C^{\infty}(\Omega\times\mathbb{R}^{n})$. Clearly, $Q$ is on the conditions of Theorem \ref{Theorem 5.1}.
\end{proof}
\begin{lemma}\label{Lemma 4.13}
Let $\Omega$ be a suitable subset of $\mathbb{H}^{n}$. Take $u:\Omega\rightarrow\mathbb{R}$ such that $\Gamma_{u}$ be a graphical translator in $\mathbb{H}^{n}\times_{-1}\mathbb{R}$. Take $\sigma$ an isometry of $\mathbb{H}^{n}$. Given $\hat{u}:=u\circ\sigma$, then $\Gamma_{\hat{u}}$ is also a graphical translator in $\mathbb{H}^{n}\times_{-1}\mathbb{R}$.
\end{lemma}
\begin{proof}
	We take $W=\sqrt{1-\lvert\nabla u \rvert ^{2}}$ such that $\mathrm{div}\left(\frac{\nabla u}{W}\right)=\frac{1}{W}$. We denote $\hat{W}=\sqrt{1-\lvert\nabla\hat{u}\rvert ^{2}}$, and then a long but straightforward computation shows $\mathrm{div}\left(\frac{\nabla\hat{u}}{\hat{W}}\right)=\frac{1}{\hat{W}}$.
\end{proof}
\begin{corollary}\label{5.2}
Take $\Omega$ a bounded open domain in $\mathbb{H}^{n}$. Assume that there exists an isometry $\sigma:\mathbb{H}^{n}\rightarrow\mathbb{H}^{n}$ such that $\sigma(\bar{\Omega})=\bar{\Omega}$. Let $u\in C^{0}\left(\bar{\Omega}\right)\cap C^{2}\left(\Omega\right)$ such that $\Gamma(p)=\left(p,u(p)\right)$ is a space-like translator, and $u\circ\sigma=u$ on $\partial\Omega$. Then, $u$ is also invariant with respect to $\sigma$, that is to say, $\Gamma$ is also invariant by $\sigma\times\, id$.
\end{corollary}
\begin{proof}
We define $\hat{u}:=u\circ\sigma\in C^{0}\left(\bar{\Omega}\right)\cap C^{2}\left(\Omega\right)$. Clearly, $\hat{u}$ is also a space-like translator, that is, $Q\hat{u}=0$. Therefore, $Qu=Q\hat{u}$. Given $x\in\partial\Omega$, $\hat{u}(x)=u\left(\sigma(x)\right)=c=u(x)$ By Theorem \ref{Theorem 5.1}, $u=\hat{u}$.
\end{proof}
\begin{remark}
There are no assumptions on the topology of $\partial\Omega$.
\end{remark}


\begin{figure}[H]
\begin{center}
	\begin{tikzpicture}[x=1.5cm,y=1.5cm]
		\draw[line width=1.pt, smooth,samples=100,] (-1.1,0)--(3,0);
		\draw[line width=1.pt, smooth, samples=100] (0.65,0) circle(0.4cm);
		\draw[line width=1.pt, smooth, samples=100] (0.2,0.5) circle(0.2cm);
		\draw[line width=1.pt, smooth, samples=100] (0.2,-0.5) circle(0.2cm);
		\clip(-1.,-1.5) rectangle (2.5,1.5);
		\draw[line width=1.pt, smooth,samples=100,domain=0.0:12.566370614359172,variable=\t] 
		plot ({cos(\t r)*(1+cos(\t r))},{sin(\t r)*(1-cos(\t r))});
	\end{tikzpicture}
	\caption{An example of $\bar{\Omega}$ which is symmetric with respect to a hyperplane.}
\end{center}
\end{figure}
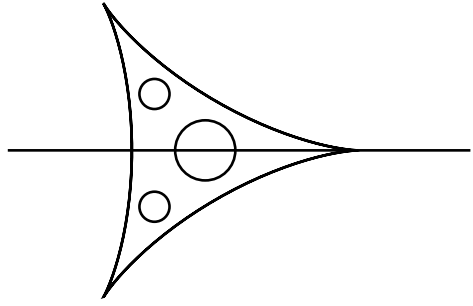
\begin{corollary}
Let $\Omega$ be a bounded open domain in $\mathbb{H}^{n}$ which is invariant by a subgroup $G$ of $\mathrm{Iso}\left(\mathbb{H}^{n}\right)$. Take $u:\bar{\Omega}\rightarrow\mathbb{R}$ such that $\Gamma_{u}$ is a space-like translator, and $u\circ\sigma=u$ on $\partial\Omega$, for all $\sigma\in G$. Then, $\Gamma_{u}$ is also invariant by $\hat{G}:=\left\{\sigma\times id:\mathbb{H}^{n}\times_{-1}\mathbb{R}\rightarrow\mathbb{H}^{n}\times_{-1}\mathbb{R}\,\vert\,\sigma\in G\right\}$.
\end{corollary}
We write $B\left(c,r\right)$ the ball of center $c$ and radius $r>0$, and $B^{*}\left(c,r\right)=B\left(c,r\right)\backslash\left\{c\right\}$.
\begin{corollary}
Let $\Omega=B^{*}\left(c,r\right)\backslash\left(c\right)\subset\mathbb{H}^{n}$, $r>0$, $u:\bar{\Omega}\rightarrow\mathbb{R}$ satisfying $\Gamma_{u}$ is a space-like translator with a possible singularity at $c$, and $u$ is constant on $\partial B\left(c,r\right)$. Then $\Gamma_{u}$ is a rotationally invariant. In addition, if $u\in C^{2}\left(B\left(c,r\right)\right)$ then, $\Gamma_{u}$ is a compact piece of the bowl.
\end{corollary}
\begin{proof}
The punctured ball $\Omega$ is symmetric with respect to all totally geodesic $\mathbb{H}^{n-1}$ passing through the center. By Theorem \ref{Theorem 5.1}, $u$ is symmetric with respect to all totally geodesic $\mathbb{H}^{n-1}$  passing through the center of $\Omega$. This means that $\Gamma_{u}$ is rotationally symmetric. So we recall Theorem \ref{4.2}. In addition, if $u$ is smooth at $c$, then $\Gamma_{u}$ has to be a piece of a bowl, bearing in mind that this is the only smooth example in the whole $\Omega$.
\end{proof}

\begin{corollary}\label{Corollary 5.4}
Take $f:\left[a_{1},b_{1}\right]\rightarrow\mathbb{R}$, $0<a_{1}<b_{1}$, one of the solutions in Theorem \ref{Theorem 3.1}. Define $c=f\left(a_{1}\right)$, $d=f\left(b_{1}\right)$ and $a_{2},b_{2}\in\mathbb{R}$, $a_{2}<b_{2}$. Then, there exists one and only one function $u:\Omega=\left[a_{1},b_{1}\right]\times\left[a_{2},b_{2}\right]\rightarrow\mathbb{R}$ such that:

1) $\Gamma_{u}$ is a space-like translator. 

2) For all $t\in\left[a_{1},b_{1}\right]$, $u\left(t,a_{2}\right)=u\left(t,b_{2}\right)=f(t)$; for all $s\in\left[a_{2},b_{2}\right]$, $u\left(a_{1},s\right)=c$, $u\left(b_{1},s\right)=d$.


\end{corollary}
Note that $\Gamma_{u}$ is foliated by horocycles as in Theorem \ref{Theorem 3.1}.
\begin{proof}
Existence is just one of the examples in Theorem \ref{Theorem 3.1}. Uniqueness: Take $u, v\in C^{0}\left(\Omega\right)\cap C^{2}(\mathring{\Omega})$ in conditions 1 and 2. We use Theorem \ref{Theorem 5.1}, so $u=v$.
\end{proof}



\section*{Declarations}

\subsection*{Funding}

Miguel Ortega is partially financed by: (1) the Spanish MICINN and ERDF, project PID2020-116126GB-I00; and (2) the “Maria de Maeztu” Excellence Unit IMAG, ref. CEX2020-001105-M, funded by MCIN/AEI/10.13039/501100011033.
and (3) Research Group FQM-324 by the Junta de Andaluc\'{i}a. \\
Buse Yalçın would like to thank the Erasmus grant from Ankara University European Union Educations Programs Coordination Office for the funding that faciliated the scientific collaboration between Ankara University (Türkiye) and University of Granada (Spain). The research of Buse Yalçın has also been supported by The Scientific and Techological Research Council of Türkiye (TÜBİTAK) Grant 2210-A.

\subsection*{Acknowledgements}
Buse Yalçın is grateful to the Institute of Mathematics (IMAG) and the Department of Geometry and Topology of Granada University for their hospitality.
\subsubsection*{Ethics declarations}

The authors declare that they have no interests/competing interests. No AI has been used to prepare this document. All pictures made by the authors of this paper with the aid of Maxima \cite{wxMaxima}.

 Miguel Ortega \\
 Department of Geometry and Topology, Faculty of Sciences, \\
 Institute of Mathematics IMAG \\
 Universidad de Granada \\
 18071 Granada \\
 Spain \\
 e-mail: miortega@ugr.es \\
 
 \noindent
 Buse Yalçın \\
 Department of Mathematics, Faculty of Sciences, \\
 Graduate School of Natural and Applied Sciences, \\
 Ankara University \\
 06100 Tandoğan, Ankara \\
 Türkiye \\
 e-mail: bsyalcin@ankara.edu.tr


\begin{thebibliography}{00}

\bibitem{AA} Alekseevsky,~A.~V., Alekseevsky,~A.~V.: Riemannian G-manifold with one-dimensional orbit space. Ann. Glob. Anal. Geom. \textbf{11}, 197–211 (1993).

\bibitem{AW}  Altschuler,~S.~J., Wu,~L.~F.: Translating surfaces of the non-parametric mean curvature flow with prescribed contact angle. Calc. Var \textbf{2} pp 101--111  (1994) https://doi.org/10.1007/BF01234317

\bibitem{BL} Batista,~M., de~Lima,~H.~F.: Spacelike translating solitons in Lorentzian product spaces: nonexistence, Calabi-Bernstein type results and examples.  Commun. Contemp. Math. \textbf{24} no. 8, 2150034, 20 pp (2022).
https://doi.org/10.1142/S0219199721500346

\bibitem{Bueno1} Bueno,~A.: Translating solitons of the mean curvature flow in the space $\mathbb{H}^2\times\mathbb{R}$. J. Geom. \textbf{109} 42, (2018).    https://doi.org/10.1007/s00022-018-0447-x

\bibitem{Bueno2} Bueno,~A.: Uniqueness of the translating bowl in $\mathbb{H}^2\times \R$. J. Geom. \textbf{111}, 43 (2020). https://doi.org/10.1007/s00022-020-00555-2

\bibitem{CSS} Clutterbuck,~J., Schnürer,~O.~C., Schulze,~F.: Stability of translating solutions to mean curvature flow. Calc. Var. \textbf{29}, 281–293 (2007). https://doi.org/10.1007/s00526-006-0033-1

\bibitem{libroEDP} Gilbarg,~D., Trudinger,~N.~S.: Elliptic partial differential equations of second order. Reprint of the 1998 edition. Classics in Mathematics. Springer-Verlag, Berlin, 2001. xiv+517 pp. ISBN: 3-540-41160-7

\bibitem{HIMW} Hoffman,~D., Ilmanen,~T., Mart\'in,~F.,  White,~B.: Notes on translating solitons for mean curvature flow. In Minimal surface: integrable systems and visualisation. vol 349 of Springer Proc. Math. Stat.. pp 147-168. Springer, Cham, 2021

\bibitem{coreano} Kim,~D.: Rotationally symmetric space-like translating solitons for the mean curvature flow in Minkowski space. J. Math. Anal. Appl. \textbf{488}, Issue 2, (2020). 124086, doi: https://doi.org/10.1016/j.jmaa.2020.124086

\bibitem{LO} Lawn,~M.-A., Ortega,~M.: Translating Solitons in a Lorentzian Setting, Submersions and Cohomogeneity One Actions. Mediterr. J. Math. \textbf{19}, 102 (2022). https://doi.org/10.1007/s00009-022-02020-7

\bibitem{LM} de~Lira,~J.~H., Mart\'in,~F.: Translating solitons in Riemannian products. J. Diff. Equations \textbf{266}, Issue 12 (2019) 7780--7812. https://doi.org/10.1016/j.jde.2018.12.015

\bibitem{MSHS15} Mart\'in,~F., Savas-Halilaj,~A., Smoczyk,~K.: On the topology of translating solitons of the mean curvature flow. Calc. Var. \textbf{54}, 2853–2882 (2015). https://doi.org/10.1007/s00526-015-0886-2

\bibitem{ON} O'Neill,~B.: Semi-Riemannian geometry, With applications to relativity. Pure and Applied Mathematics, 103. Academic Press, Inc. New York, 1983.


\bibitem{Sol} Pipoli,~G. Invariant translators of the solvable group. Annali di Matematica \textbf{199}, 1961–1978 (2020). https://doi.org/10.1007/s10231-020-00951-0

\bibitem{Heis} Pipoli,~G. Invariant Translators of the Heisenberg Group. J Geom Anal \textbf{31}, 5219–5258 (2021). https://doi.org/10.1007/s12220-020-00476-1

\bibitem{RS} Ros,~A., Sicbaldi,~P.: Geometry and topology of some overdetermined elliptic problems. J. Differential Equations \textbf{255}(2013), no.5, 951–977. https://doi.org/10.1016/j.jde.2013.04.027

%

\bibitem{Wiggins} Wiggins,~S.:
Introduction to applied nonlinear dynamical systems and chaos.
Second edition. Texts in Applied Mathematics, 2. Springer-Verlag, New York, 2003. ISBN: 0-387-00177-8

\bibitem{wxMaxima} wxMaxima,  \href{https://maxima.sourceforge.io/}{https://maxima.sourceforge.io/}
Last accessed: 2023-July-13.


\end{thebibliography}
\end{document}